\documentclass[12pt]{article}
\usepackage{amsfonts}

\begin{document}
\def\C{{\mathbb{C}}}
\def\R{{\mathbb{R}}}
\def\H{{\mathbb{H}^n}}

\setlength{\parindent}{8mm}
\setlength{\oddsidemargin}{0pt}
\setlength{\evensidemargin}{0pt}
\def\grad {\nabla}
\def\theequation{\thesection.\arabic{equation}}
\def\SQ{\displaystyle{\mathrel{\mathop{S_2^Q}}}}
\def\N{{\rm I}\!{\rm N}}
\def\G{{\b G}}
\def\T{{\cal T}}
\def\Dh{\Delta_{\H}}
\def\Dg{\Delta_{\H}}
\def\ur{u_{\cal  S}}
\def\grad{\nabla}
\def\p{\partial}
\def\la{\lambda}
\def\l{\lambda}
\def\S{\displaystyle{\mathrel{\mathop{S_1^2}^{o}}}}
\newtheorem{thm}{Theorem}[section]
\newtheorem{pro}{Proposition}[section]
\newtheorem{lem}{Lemma}[section]
\newtheorem{cor}{Corollary}[section]
\newtheorem{defi}{Definition}[section]
\title{A negative answer to a one-dimensional symmetry problem in the
Heisenberg group}
\author{I. Birindelli\thanks{Dipartimento di
Matematica,
Universit\`a di Roma "La Sapienza",Piazzale A.Moro 2 - I 00185 Roma,
Italy;
  e-mail: isabeau@mat.uniroma1.it
}, E. Lanconelli\thanks{Dipartimento di Matematica, Universit\`a di
Bologna, Piazza di Porta San Donato 5, I 40126 Bologna;
 e-mail: lanconelli@mat.unibo.it}}
\date{}
\maketitle

\section{Introduction and Main Results}
Symmetry properties of solutions to semi-linear elliptic equations have
been widely studied in the last decades. In this contest, a long
standing conjecture by De Giorgi states that any global solution
to the Ginzburg-Landau equation
 $$
 \Delta u + u(1-u^2)=0 \
\mbox{in}\ \R^N \ (N\leq 8)
 $$
satisfying $-1\leq u\leq 1$ and $\frac{\partial u}{\partial
x_N}>0$ is constant along hyperplanes. Recently this conjecture
was proved to be true by Ghoussoub and Guy for $N=2$ (\cite{gg})
and by Ambrosio and Cabr\'e for $N=3$ (\cite{ac}). It is still an
open question for $N>3$ though Alberti, Ambrosio and Cabr\'e
generalized the result for any $C^2$ non-linearity (when $N\leq
3$) \cite{aac}.

Under the further hypothesis that the solution $u$ satisfies
 $$
\lim_{x_3\rightarrow \pm \infty} u(x',x_3)=\pm 1\ \ \forall x'\in
\R^2
  $$
the proof that $u$ is constant along hyperplanes given in
\cite{ac} is somehow simpler. On the other hand, under the
hypothesis that this limit is uniform in $x'$, the conjecture was
known as Gibbons conjecture and it has been proved for all
dimensions independently by Barlow, Bass, Guy in \cite{bbg},
Berestycki, Hamel , Monneau in \cite{bhm} and Farina in \cite{f}.

In recent years symmetry and monotonicity properties of solutions
to semilinear equations have been investigated in the more general
contest of the Carnot groups, see \cite{bp1, bp2,bp3}, \cite{ag},
\cite{bl} and \cite{gv}. The interest in semi-linear equations in
these groups has increased as they appear in many theoretical and
application fields, such as complex geometry and mathematical
models for crystal structures \cite{c}.

In \cite{bp2} Prajapat and the first author studied Gibbons
conjecture for the equation
 \begin{equation}\label{1.3}
 \Dg u+f(u)=0 \ \mbox{in} \  \H,
 \end{equation}
where $\Dg$ denotes the Kohn-Lalacian on the Heisenberg group $\H$
and $f(u)$ is a non linear term with some general hypothesis (in
particular they include the case $f(u)=u(1-u^2)$). They prove that
the conjecture holds true for all directions orthogonal to the
center of $\H$. \footnote{Very recently, in \cite{bl}, the results
of \cite{bp2} have been extended to every sub-Laplacian on a
Carnot group.} The question of whether the result holds true in
the remaining direction was raised in \cite{bp2}.

The aim of this paper is to prove that, with respect to the center
direction of $\H$, the stronger De Giorgi conjecture is not true
for the equation (\ref{1.3}). This negative answer will easily
follow from next Theorem 1.1, the main result  of this note.

In order to clearly state our theorem, we need to recall some
known facts about the Heisenberg space $\H$ and its intrinsic
Laplacian $\Dg$.

First of all let us say that  $\H$ is the Lie group whose
underlying manifold is $\C^n\times \R$, $n\in N$, endowed with the
group action $\circ$ given by
\begin{equation}
 \label{ga}
 \xi_o\circ \xi= (z+z_o,t+t_o +2{\rm Im}(\overline{z}\cdot z_o) ).
\end{equation}
Here and in the rest of the paper we identify $\C^n$ with
$\R^{2n}$ and, setting $z=x+iy$, for the point of $\H$ we use
the equivalent notations 
$\xi=(z,t)=(x,y,t)\in \R^n\times\R^n\times\R$ with
$z:=(z_1,\dots,z_n)=(x_1,y_1,\dots,x_n,y_n)$. Furthermore, ``$\,\cdot\,$''
denotes the usual inner product in $\C^n$.

The Lie Algebra of left invariant vector fields is generated by
\[ \begin{array}{l}
X_i=\frac{\partial}{\partial x_i} + 2y_i \frac{\partial}{\partial
t}, \;\;\mbox{for }\;\;i=1,\dots,n,\\ Y_i=\frac{\partial}{\partial
y_i} - 2x_i \frac{\partial}{\partial t},\;\;\mbox{for
}\;\;i=1,\dots,n.
\end{array}
\]
The intrinsic Laplacian of $\H$, also called the Kohn Laplacian, is
defined as
 $$\Dg=\sum_{i=1}^{n} (X_i^2+Y_i^2).$$
 It is a second order
degenerate elliptic operator of Hormander type and hence it is
hypoelliptic (see e.g. \cite{fs} or \cite{j} for more details
about $\Dg$).

 With respect to the
group dilation $\delta_{\l}\xi=(\l z,\l^2t)$,  $\Dg$ is
homogeneous of degree two in the following sense
\[ \Dg\circ \delta_\l=\lambda^2 \delta_\l\circ\Dg.\]

The Koranyi ball of center $\xi_o$ and radius $R$ is defined by
\[B_H(\xi_o,R):=\{\xi\;\mbox{such that } |\xi^{-1}\circ\xi_o|_\H\leq R\}\]
where
 $$|\xi|_{\H}=\left( |z|^4 + t^2\right)^{\frac{1}{4}}$$ is a
norm with respect to the group dilation and it satisfies
\[ |B_H(\xi_o,R)|=|B_H(0,R)|=CR^Q\]
where $Q=2n+2$ is the homogeneous dimension of $\H$.

A fundamental solution of $-\Dg$ with pole at the origin is given
by:
 $$
 \Gamma(\xi)=\frac{C_Q}{(|\xi|_{\H})^{Q-2}}
$$
 where $C_Q$ is a positive constant.

\medskip
For our purposes it is convenient to remind that the class of
cylindrically symmetric functions is  invariant with respect to
the action of $\Dg$. We shall say that a function
$(z,t)\rightarrow u(z,t)$ is cylindrically symmetric if there
exists a two variables function $U$ such that $u(z,t)=U(r,t)$,
$r=|z|$.

In that case we formally have that
 $$ \Dg u(z,t)= \partial_{rr}U+\frac{2n-1}{r}\partial_r U
 +4r^2\partial_{tt}U.$$
The main result of this paper is the following:
\begin{thm} Let $f:\R\rightarrow \R$ be a locally Lipschitz
continuous function satisfying the hypotheses listed below:

(H1) $f$ is odd,

(H2) $f>0 $ in $]0,1[$, $\, f(0)=f(1)=0$,

(H3) $\displaystyle\lim_{s\rightarrow 0} \frac{f(s)}{s}=l>0$.

Then there exists a solution $u$ to the equation:
 \begin{equation}\label{1.4}
 \Dg u+f(u)=0 \quad \mbox{in}\quad \R^{2n+1}
 \end{equation}
 satisfying $|u|< 1$, $ \frac{\partial u}{\partial
t}>0$ and
 $$ \displaystyle\lim_{t\rightarrow \pm \infty}
u(z,t)=\pm 1.
 $$
  Moreover $u$ is cylindrically symmetric and of
class $C^\infty$ when $f$ is $C^\infty$.
\end{thm}
For solution $u$ of (\ref{1.4}) we mean a continuous function $u$
such that:
\begin{enumerate}
\item For a suitable $\alpha>0$, $u\in\Lambda^{2+\alpha}_{\rm
 loc}(\H)$ i.e. $X_j^2u$ and $Y_j^2u$, $j=1,\cdots,n$, exists in
 the weak sense of distributions and belong to $\Lambda^{\alpha}_{\rm
 loc}(\H)$

\item $u$ satisfies  (\ref{1.4}) pointwise everywhere.
\end{enumerate}

As in \cite{fs} we have denoted by $\Lambda^{\alpha}_{\rm loc}(\H)$
the class of functions which are locally $\alpha$-Holder
continuous with respect to the intrinsic distance $d$ in $\H$
defined by
 $$d(\xi,\xi')=|(\xi')^{-1}\circ \xi|_\H.$$
Using the commutators of the Lie Algebra,
 it is easy to see  that
$\Lambda^{2+\alpha}_{\rm loc}(\H)$ is continuously embedded in the
usual $C^{1+\frac{\alpha}{2}}_{\rm loc}(\R^{2n+1})$.

\medskip
 From Theorem 1.1  we immediately get the following corollary.

\begin{cor} De Giorgi's conjecture in the $t$-direction is not true in
$\H$.
\end{cor}
{\bf Proof.} The functions $f(s)=s(1-s^2)$ satisfies all
hypotheses of Theorem~1.1, hence there exists a $C^\infty$
function $u$ such that
 $$
 \left\{
 \begin{array}{lc}
 \Dg u +u(1-u^2)=0 & \mbox{in }\ \R^{2n+1},\\
 -1<u<1,  \, \frac{\partial u }{\partial t}>0, &\\
 \displaystyle\lim_{t\rightarrow \pm \infty}u(z,t)=\pm 1.
 \end{array}
 \right.
 $$
Then, if De Giorgi conjecture were true in the $t$ direction there
would exist $\alpha\in \R^{2n}$ and $\nu>0$ such that
$u(z,t)=U(\alpha \cdot z +t\nu)$ for some function
$U:\R\rightarrow \R$. Furthermore $U$ would satisfy
 $$\left(|\alpha|^2 - 4\nu (J\alpha\cdot z)
 +4r^2\nu^2\right)U''=U(U^2-1)$$
where $J$ is the classical symplectic $2n\times 2n$ matrix. This
is a contradiction since the right hand side is constant along the
hyperplanes $\alpha \cdot z +t\nu=c$ for any $c\in \R$ while the
left hand side is not.

\bigskip
\noindent{\bf Remark 1.1} It would be interesting to know whether the
function constructed in Theorem 1.1 has uniform limit with respect
to $z$.

\noindent{\bf Remark 1.2} It is natural to consider the extension
of Theorem 1.1 to the contest of Carnot groups. This will be the
object of a subsequent study.

\section{Proof of Theorem 1.1.}
For any $R>0$ we shall denote by $D_R$ and $D_R^+$ respectively
the cylinders
 $$
D_R=\{(z,t)\in \R^{2n+1} ; \ |z|<R, \ |t|<R^2\} $$ and
 $$
D_R^+=\{(z,t)\in \R^{2n+1} ;\  |z|<R, \ 0<t<R^2\}. $$ Let
$\psi(t)=\frac{t}{R^2}$.

We shall split the proof in several steps.

\noindent {\bf First step:} {\em The semilinear Dirichlet problem
\begin{equation}\label{i2}
\left\{\begin{array}{l} \Dg u = -f(u) \mbox{~~in~~}D_R^+ ,\\
u(r,t)=\psi(t),\; \mbox{~~on~~} \partial D_R^+.
\end{array}
\right.
\end{equation}
has a solution $u\in \Lambda^{2+\alpha}_{\rm loc}(D_R^+)\cap
\Lambda^{\alpha}(\overline{D_R^+})$ for a suitable $\alpha\in
(0,1)$. Furthermore $u$ is cylindrically symmetric,  $0\leq u\leq
1$ and for any $R$ sufficiently large,
 $$u\geq v_o$$
for some function $v_o\geq 0$, $v_o\not\equiv 0$, $v_o$
independent of $R$.}

\bigskip
Let $M\in\R^+$ be larger than the Lipschitz constant of $f$ in
$[0,1]$ and let us define
 $$g:\R\rightarrow \R, \ g(s)=f(s)+Ms.$$


Let $\T$ be the map formally defined by
$\T(v)=u$ where $u$ is the only solution to the Dirichlet problem

\begin{equation}\label{i3}
\left\{\begin{array}{l} \Dg u -Mu= -g(v) \mbox{~~in~~}D_R^+ ,\\
u=\psi,\; \mbox{~~on~~} \partial D_R^+.
\end{array}
\right.
\end{equation}

The operator $\T$ has the following properties:

\bigskip
\noindent (P1) There exists $\alpha\in (0,1)$ such that $\T$ is
well defined in $\Lambda^{\alpha}(\overline{D_R^+})$. Furthermore
 \begin{equation}\label{2.4}
 |u(\xi)-u(\xi')|\leq C d(\xi',\xi)^\alpha (1+\sup|g(v)|)
 \end{equation}
  for any $\xi, \xi'\in D_R^+$. We also have that $\T(v)\in \Lambda^{2+\alpha}_{\rm loc}
  (D_R^+)$ for every $v\in \Lambda^{\alpha}(\overline {D_R^+})$.

This statement can be proved by using standard arguments and the
results in \cite{fs}, \cite{j} (see also \cite[Theorem 4.1]{gl}).

\bigskip
\noindent (P2) $\T(v)$ is cylindrically symmetric if $v$ is
cylindrically symmetric.

Indeed suppose that $u=\T(v)$. Let $\cal S $ be a rotation in
$\R^{2n}$ and define $\ur(z,t):=u({\cal S} z,t)$. Since $\Dg$
is invariant with respect to $\cal S$,  we have $\Dg \ur(z,t)=\Dg u ({\cal
S} z,t)$, so that $\ur$ is a solution of

$$ \left\{\begin{array}{l} \Dg \ur -M\ur= -g(v({\cal S}z,t))=-g(v)
\mbox{~~in~~}D_R^+ ,\\ \ur=\psi,\;  \mbox{~~on~~} \partial D_R^+.
\end{array}
\right.
$$

Here we have used the invariance with respect to $\cal S$ of $v$,
$\psi$ and $D_R^+$ .

By the maximum principle we know that the solution of (\ref{i3})
is unique, hence $u=\ur$ for any ${\cal S}$, i.e. $u$ is
cylindrically symmetric.

\bigskip
\noindent (P3) $\T$ is monotone. More precisely if $v_1,v_2\in
\Lambda^{\alpha}(\overline {D_R^+})$ and  $0\leq v_1\leq  v_2\le 1$,
then $\T v_1\leq \T v_2$.

Let us observe that with our choice of $M$ if $0\leq v_1\leq  v_2$
then $g(v_1)\leq g(v_2)$. Hence (P3) follows from the maximum
principle for $-\Dg+ M$ in $D_R^+$.

\bigskip
\noindent (P4) If $v\in \Lambda^{\alpha}(\overline {D_R^+})$ and
$0\leq v\leq 1$ then $0\leq \T(v)\leq 1$.

Indeed, since $g(0)=0$, $g(1)=M$ and $0\leq \psi\leq 1$ on
$\partial D_R^+$, again by the maximum principle we obtain that
$\T(1)\leq 1$ and $\T(0)\geq 0$. Now we only need to apply
property (P3) for $v\in \Lambda^{\alpha}(\overline {D_R^+})$ such
that $0\leq v\leq 1$.


\bigskip
We shall now construct a function $v_o\geq 0$ that plays the role
of a lower barrier.

Let $\lambda_o$ denote the principal eigenvalue of $-\Dg$ in
$D_R^+$ and let $\phi_o>0$ be the corresponding eigenfuntion
normalized by $\sup\phi_o=1$.

We choose and  fix $R_o$ sufficiently large that
 $$\lambda_o\leq  \frac{l}{2}$$
where $l$ is the limit in condition (H3). Then  there exists
$\varepsilon\in (0,1)$ independent of $R$ such that
 $$\lambda_o\varepsilon\phi_o\leq f(\varepsilon\phi_o).$$

By uniqueness of the normalized eigenfunction $\phi_o$, arguing as
in the proof of (P2) we can prove that  $\phi_o$ is cylindrically
symmetric.

 From now on we assume that $R>R_o$. Let us define
 $$v_o=\left\{\begin{array}{lc}
             \varepsilon \phi_o & \quad\mbox{in}\quad D_{R_o}^+\\
            0  & \quad\mbox{in}\quad D_R^+\setminus D_{R_o}^+.
\end{array}
\right.$$

Standard arguments show that $v_o$ is locally Holder continuous in
$\R^{2n+1}$, (see e.g. \cite[Theorem 4.1]{gl},  we stress that
condition (4.4) in that theorem is satisfied since   $D_R^+$ is
convex).

As a consequence $\T(v_o)$ is well defined and since $0\leq
v_o\leq 1$ using (P4) we get that $0\leq \T(v_o)\leq 1$. Let us
now prove that $v_o\leq u_o:=\T(v_o)$. Clearly the inequality
holds in $D_R^+\setminus D_{R_o}^+$, using (P4), hence we just
have to prove it in $D_{R_o}^+$. We have

\begin{eqnarray}
\Dg u_o-Mu_o &= &-g(v_o)=-g(\varepsilon\phi_o)\leq
-(M+\lambda_o)(\varepsilon\phi_o)\nonumber\\
 &=& -M\varepsilon\phi_o +\Dg\varepsilon\phi_o =-Mv_o+\Dg v_o,
 \nonumber
 \end{eqnarray}
so that
$$\left\{\begin{array}{lc}
          \Dg (u_o-v_o)-M(u_o-v_o)\leq 0& \mbox{in}\  D_{R_o}^+\\
              u_o\geq v_o &\mbox{on}\ \partial D_{R_o}^+.
\end{array}
\right. $$
 The maximum principle implies that $u_o\geq v_o$ in $D_{R_o}^+$.

\bigskip

Now we construct the sequence of functions
 $$v_k=\T^k(v_o), k\in\N. $$

Clearly using the properties above, all $v_k$ are cylindrically
symmetric  and
 $$1\geq \T^k(v_o)\geq \T(v_o)\geq v_o\geq 0\ \mbox{for every }\ k\in\N.$$
Let us denote by $u$ the pointwise limit of  $(v_k)$. Then $u$ is
cylindrically symmetric , $v_o\leq u\leq 1$, $u\in
\Lambda^{\alpha}(\overline {D_R^+})$ since, by (\ref{2.4})

$$ |v_k(\xi)-v_k(\xi')|\leq C d(\xi',\xi)^\alpha $$
 where $C>0$ is independent of $R$. This estimates implies that
 the $v_k$ uniformly converges to $u$ in $\overline{ D_R^+}$, so
 that $u=\psi$ on $\partial D_R^+$.

Furthermore in the weak sense of distributions, $u$ satisfies
\begin{equation}\label{2.6}
\begin{array}{l} \Dg u +f(u)=0
\mbox{~~in~~}D_R^+ .
\end{array}
\end{equation}
 From (\ref{2.6}), the Holder regularity of $u$ and standard
bootstrap argument we obtain that $u\in \Lambda^{2+\alpha}_{\rm
loc} (D_R^+)$ and it satisfies the equation pointwise. Hence $u$
is the required function.

\noindent{\bf Remark 2.1} Since $u$ is cylindrically symmetric we
have that $u(z,t)=U(|z|,t)$ and $U$ satisfies the semilinear
elliptic equation
 $$\partial_{rr}U+\frac{2n-1}{r}\partial_{r}U +4r^2 \partial_{tt}U
 +f(U)=0$$
 in the open subset of $\R^2$
 $$ \Omega_R:=\{(r,t)\in \R^2/ \ 0<r<R, 0<t<R^2\}.$$
Moreover $U$ is locally $\frac{\alpha}{2}$-Holder continuous, in the usual
sense, up to $\partial \Omega_R\setminus \{(0,t)/ \ 0<t<R^2\}$.
Then, being $U(r,0)=0$ when $0<r<R$, by classical regularity
results for elliptic equations, $U$ is of class
$C^{2+\frac{\alpha}{2}}_{\rm loc}$ up to $\Omega_R\cup \{(r,0)/\
0<r<R\}.$

\bigskip
\noindent {\bf Second step:} {\em The function constructed in the
first step satisfies} $\frac{\partial u}{\partial t}> 0$.

\bigskip
In \cite{bp3} the following definition  and theorem are given:
\begin{defi} Fix $\eta \in \H$.
A domain $\Omega \subset H$ is said to be $\eta$-convex (or convex in
the direction $\eta$) if for any
$\xi_1 \in \Omega$ and any $\xi_2 \in \Omega$ such that $\xi_2 = \alpha \eta \circ \xi_1$ for some
$ \alpha> 0 $, we have $s  \eta \circ \xi_1 \in \Omega$ for every $s \in
(0,\alpha)$.
\end{defi}

\begin{thm}\label{mono}
Let $\Omega$
be an arbitrary bounded domain of $\H$ which is $\eta$- convex for
some $\eta \in H$. Let $ u \in \SQ (\Omega) \cap C(\bar{\Omega})$
be a solution of
\begin{equation}\label{main}
\left. \begin{array}{rll}
\Dg u + f(u) & = & 0 {\rm ~~ in ~~} \Omega\\
        u & = & \phi{\rm ~~ on~~} \partial \Omega
\end{array}\right\}
\end{equation}
where $f$ is a Lipschitz continuous function. Assume that for any
$\xi_1$, $\xi_2 \in \p \Omega$, such that $ \xi_2 = \alpha \eta \circ
\xi_1$ for some $\alpha> 0$,  we have for each $s \in (0, \alpha)$
\begin{equation}\label{1}
\phi(\xi_1) < u(s\eta\circ {\xi_1}) < 
\phi(\xi_2){\ ~~} s \eta \circ \xi_1 \in  \Omega
\end{equation}
and
\begin{equation}\label{2}
\phi(\xi_1) < \phi(s\eta {\xi_1}) < \phi(\xi_2){\rm~~ if ~~}
  s\eta \circ \xi_1 \in  \partial\Omega
\end{equation}

Then $u$ satisfies

\begin{equation}{\label{result}}
u (s_1\eta\circ {\xi}) < u (s \eta\circ {\xi})
 \end{equation}
for any $0<s_1 < s<\alpha$ and  for every $\xi \in \Omega$.

Moreover, $u$ is the unique solution of (\ref{main})
 in  $\SQ(\Omega) \cap C(\bar{\Omega})$ satisfying (\ref{1}).
\end{thm}

Let us choose $\eta=(0,1)$,  clearly $D_R^+$ is $\eta$-convex
since: $$s \eta\circ \xi=(z,t+s).$$

Furthermore $0=\psi(0)\leq u(z,t)\leq \psi(1)=1$ and by
construction $\psi$ satisfies (\ref{2}). Hence we are in the
hypothesis of Theorem \ref{mono} and $u$  satisfies

$$u(z,t_1)\leq u(z,t_2)\ \mbox{ for any } 0\leq t_1\leq t_2\leq 1$$
 in $D_R^+$.

In particular we get $\frac{\partial u}{\partial t}\geq 0$.

Now since $\frac{\partial }{\partial t}$ commutes with $\Dg$ and
$f$ is Lipschitz continuous then the inequality is strict, just by
using the strong Maximum principle.

\bigskip
\noindent {\bf Third step.} {\em We extend to $D_R$ the function
$u$ of the previous step by setting}
$$v(z,t)=\left\{\begin{array}{lc}
                   u(z,t) & \mbox{ for }\ t\geq 0\\
                   -u(z,-t)     & \mbox{ for }\ t\leq 0.
                    \end{array}\right.
$$

Obviously $v$ is cylindrically symmetric, $-1\leq v\leq 1$, $v\geq
v_o$ in $D_R^+$, $v\in C^{\frac{\alpha}{2}}(D_R)$ and $v=\phi$ on
$\partial D_R$. We want to prove that $v$ satisfies
 \begin{equation}\Dg v+f(v)=0 \ \mbox{in }\
 D_R.\label{2.7}\end{equation}
Since  $f$ is odd, using the fact that $v$ is odd and
cylindrically symmetric it is easy to see that  $v$ satisfies
(\ref{2.7}) in $D_R\setminus\{t=0\}$.

By Remark 2.1 at the end of the first step, we now obtain that
$v\in C^{2+\frac{\alpha}{2}}(D_R\setminus \{(0,0)\} )$ and it
solves (\ref{2.7}) in the same open set. Hence we just have to
remove the singularity at the origin. Let us define
 $$w(\xi)=-\int_{D_R}\Gamma((\xi')^{-1}\circ\xi) f(v(\xi'))d\xi',$$
 where
$\Gamma(z,t)$ is the fundamental solution recalled in the
Introduction. Since $f(v)\in C^{\frac{\alpha}{2}}( D_R )$
and
 $C_{loc}^{\frac{\alpha}{2}}(D_R )\subset\Lambda_{loc}^{\frac{\alpha}{2}}(D_R)$, then $w\in \Lambda_{loc}^{2+\frac{\alpha}{2}}(D_R)$ and satisfies
 $$\Dg w=f(v) \ \mbox{in} \ D_R.$$
Hence $$\Dg (v+w)=0 \ \mbox{ in } \ D_R\setminus\{(0,0)\}.$$

On the other hand $v+w\in L^{\infty}(D_R)$. Then there exists a
$C^{\infty}$-function $h$, $\Dg$ harmonic in $D_R$ such that
 $$h=v+w\ {\rm in} \ D_R\setminus \{(0,0)\}.$$

It follows that $v$ solves (\ref{2.7}) everywhere in $D_R$.

This ends the third step. We shall denote $u_R(z,t)=v(z,t)$ the
function constructed above.

\bigskip
\noindent{\bf Fourth step.} {\em We let $R$ tend to infinity and
obtain a global solution.}

Since the functions  $u_R$ are equi-bounded and solutions of
(\ref{i2}) in $D_{R}$, then $\Dg u_R$ are also equi-bounded and by
standard arguments, eventually passing to a subsequence, the $u_R$'s
locally uniformly converge to $u$, weak solution of
\begin{equation}\label{rn}\Dg u+f(u)=0  \  \mbox{in }\ \R^{2n+1}.
\end{equation}

\noindent Furthermore

1)$u$ is cylindrically symmetric,

2) $-1\leq u\leq 1$,

3) $u(z,t)=-u(z,-t)$,

4) for $t\geq 0$, $u(z,t)\geq v_o(z,t)$,

5) $t\mapsto u(z,t)$ is monotone increasing.

\noindent
Since $f$ is locally Lipschitz continuous and $|u|\le 1$, it follows
from (\ref{rn}) that $u\in \Lambda_{loc}^{2+\alpha}(\H)$ for
every $\alpha<1$. Obviously, the more regular $f$ is, the more
regular $u$ is; in particular $u$ is of class $C^\infty$ when
$f$ is $C^\infty$. 

Moreover, property  5) implies $\frac{\partial u}{\partial t}\geq 0$
so that, since $\frac{\partial}{\partial t}$ commutes with $\Dg$,
by the strong maximum principle either $\frac{\partial u}{\partial
t}>0$ or $\frac{\partial u}{\partial t}\equiv 0$. But by 3) and 4)
this second possibility is absurd hence $\frac{\partial
u}{\partial t}>0$ .

\bigskip
\noindent{\bf Last  step.} {\em We want to prove that}
$$\lim_{t\rightarrow \pm \infty} u(z,t) =\pm 1.$$

We shall consider only the limit in $+\infty$ since the other case
follows similarly. Let us denote
$u_o(z):=\displaystyle\lim_{t\rightarrow +\infty} u(z,t)$. Since
$u$ is bounded and monotone in $t$ the limit is well defined and
$0<u_o(z)\leq 1$. We want to prove that $u_o(z)\equiv 1$.

By standard arguments (multiplying equation (\ref{rn}) by a
sequence of functions $\psi_k(z,t)=\phi(z)\phi_k(t)$ where $\phi$
has compact support and supp$\phi_k =]k,k+1[$ and $\int \phi_k
dt=1$ and letting $k$ go to infinity) it easy to see that $u_o$ is
a weak solution of

$$\Delta u_o+f(u_o)=0 \ \mbox{ in}\  \R^{2n}.$$

Clearly a bootstrap argument shows that $u_o$ is a classical
solution. Moreover $u_o(z)=U_o(r)$ with $r=|z|$  for some function
$U_o$ solution of
\begin{eqnarray}\label{er}
U_o''(r)+\frac{2n-1}{r}U_o'(r) + f(U_o(r)) &=& 0, \\
 U_o'(0)=0 &&
\end{eqnarray}
The Cauchy problem for (\ref{er}) with initial conditions
$U_o(0)=1$ and $U_o'(0)=0$ has a unique solution (see e.g.
\cite{ps}). Thus,  since $f(1)=0$, if $U_o(0)=1$ then $U_o\equiv 1$
and we are done. Suppose, by contradiction, that $U_o(0)<1$.

It is easy to see that $U_o'<0$.
 Indeed integrating (\ref{er}) one obtains:
\begin{equation}\label{er2} 
r^{2n-1}U_o'(r)=-\int_0^r \rho^{2n-1}f(U_o(\rho))d\rho<0,
\end{equation}
hence $U_o$ is strictly decreasing and has a finite
non--negative limit as $r\rightarrow~\infty$. More precisely
$\displaystyle\lim_{r\rightarrow +\infty}U_o(r)=0$. Indeed
otherwise $U_o(r)\rightarrow k>0$
 and $f(U_o(r))\rightarrow f(k)>0$ (by (H3)). This, together with
(\ref{er2}) implies
that $|U_o'(r)|\rightarrow \infty$, which is absurd since $U_o$ is bounded.
Using hypothesis (H4) on $f$ we obtain that for $r$ large $U_o$
satisfies $$U_o''(r)+\frac{2n-1}{r}U_o'(r)+K(r)U_o(r)=0$$ with
$K(r)=\frac{f(U_o(r))}{U_o(r)}\rightarrow l>0$.

Using the substitution $V_o(r)=r^{\frac{2n-1}{2}}U_o(r)$ we obtain
that $V_o$ satisfies

$$V''(r)+H(r)V(r)=0$$ with
$H(r)=\frac{2n-1}{2}(1-\frac{N-1}{2})\frac{1}{r^2} +K(r)$.
Comparing with $$U''(r)+\frac{l}{2}U(r)=0$$ we obtain that $V_o$
i.e. $U_o$ has infinite zeros in a neighborhood of infinity, which
is absurd. This conclude the last step and the proof.

\end{document}